\documentclass[12pt,a4paper]{article}
\usepackage{amsmath,amssymb,amsfonts, amsthm}
\date{\begin{flushleft} {\small } \end{flushleft}}
\usepackage[cp1251]{inputenc}
\usepackage[english,russian]{babel}
\usepackage{amssymb,amsmath}

\newtheorem{thm}{Theorem}[section]
\newtheorem{ut}{Statement}
\newtheorem{utt}{Proposition}[section]
\newtheorem{cor}{Corollary}[section]
\newtheorem{lm}{Lemma}[section]
\theoremstyle{remark}
\newtheorem{rem}{Remark}[section]
\theoremstyle{definition}
\newtheorem{df}{Definition}[section]

\newtheorem{ex}{Example}[section]
\newtheorem{qu}{Question}

\numberwithin{equation}{section}

 \def\proof{{\textbf{Proof}. }}
 \def\T{{\mathbb T}}
\def\C{{\mathbb C}}

\def\R{{\mathbb R}}

\def\h{{\rm h}}
\def\re{{\rm Re}}
\def\im{{\rm Im}}

\def\codim{{\rm codim\:}}
\def\codimc{{\rm codim}_\C\:}
\def\supp{{\rm supp\:}}

\def\min{{\rm min}}

\def\1{{1}_\T}

\def\({{\rm (}}
\def\){{\rm )}}
\def\st{{\rm st}}

\def\spr{{\cap^\st}}

\def\ETP{{{\rm ETP}}}
\def\ETV{{{\rm ETV}}}

\def\Int{{\rm Int}}

\def\mdskip{\vskip-\lastskip\vskip\medskipamount}

\def\QED{{\parfillskip0pt\hfil$\square$\par}\mdskip}
\author{B. Kazarnovskii}
\title{Exponential tropical varieties and complex Monge-Ampere operator\thanks {Supported by NSh-4850.2012.1}}
\begin{document}
\maketitle {\small{
Sometimes it is possible to extend some using Newton polyhedra computations in algebraic geometry from polynomials
to exponential sums.
For this purpose it is useful to consider analogues of tropical varieties
in complex space.
These analogues are called exponential tropical varieties (ETV).
We construct the ring of ETV.
Algebraic tropical varieties form the subring of the ring of ETV.
In this paper we connect ETV with the complex Monge-Ampere operator action on the space of piecewise
linear functions in complex vector space.
We show that all ETV arise as results of such operator action.
We give some applications of this connection.
One of the applications is a criterion for zero value of a mixed Monge-Ampere operator.
This criterion is the modification of the criterion for zero value of a mixed volume of convex bodies.
The proof is the modification of A. Khovanskii's unpublished proof of the corresponding theorem on mixed volumes.
In the part 1 we give the definition of ETV and detail statements of theorems (without proofs).
In the part 2 we prove the theorems on the action of the Monge-Ampere operator.
}
}
\selectlanguage{english}
\tableofcontents
\section{Definition and basic properties of ETV}
%
%
 \subsection{The definition of ETV}
  Let $X$ be a finite set of closed convex
   (not necessarily bounded) $k$-dimensional polytopes
   in a real vector space $E$.
  We say that $X$ is a $k$-dimensional polyhedral set,
   if intersection of any two polytopes
   either is empty or is their common face.
   Any face of any polytope is called a \emph{cell}.
  By default $E$ is the space $\C^n$  considered as a real vector space.
  Also, we identify the dual space $E^*$ with ${\C^n}^*$
  using the pairing $(z,z^*)=\re\langle z,z^*\rangle$.
  By definition, \emph{a chain of degree} $m$ on $X$ is the odd function
   on the set of oriented $k$-dimensional cells taking $\Delta$ to $X_\Delta\in\bigwedge_\C^m {\C^n}^*$
   (the function is called odd, if it's value at the argument $\Delta$ changes sign with the changing of orientation of $\Delta$).
  The polyhedral set with a fixed chain we call \emph{a framed polyhedral set}.
  The odd form $X_\Delta$ we call the frame of a cell $\Delta$.
  \begin{df}
 The union of $k$-dimensional cells $\Delta\in X$ with nonzero frames $X_\Delta$
 is denoted $\supp X$ and is called \emph{the support} of the framed polyhedral set.
 Say that two framed $k$-dimensional polyhedral sets $X,Y$ with the common support
\emph{are equivalent},
 if $X_\Delta=Y_\Theta$
 for any $k$-dimensional cells $\Lambda\in X,\Theta\in Y$
 with $k$-dimensional intersection.
  \label{dfEquPol}
 \end{df}
   The set of $(k-1)$-dimensional cells of framed $k$-dimensional polyhedral set $X$ form the $(k-1)$-dimensional polyhedral set
   $\partial X$.
   Make it framed as
    \begin{equation}
    \label{bnd}
     (\partial X)_\Lambda=\sum_{\Delta\supset\Lambda,\dim\Delta=k} X_\Delta,
   \end{equation}
    where the orientations of the cells $\Lambda$ and $\Lambda$ agreed as usual.
    The framed polyhedral set is called closed if $\supp \partial X=\emptyset$.
 \begin{cor}
  The framed polyhedral set $\partial X$ is closed.
 \label{corPolBound}
 \end{cor}
  Let $E_\Delta$ be a tangent space of the cell $\Delta$
  and $\C_\Delta$ be a maximal complex subspace of $E_\Delta$.
  \begin{df}
 Let $k\geq n$ and $X$ be a closed framed $k$-dimensional polyhedral set with the chain of degree $2n-k$.
 Say that $X$ is an exponential tropical polyhedral set (ETP),
 if for any $k$-dimensional cell $\Delta$

 {\rm(1)} the restriction $X_{\Delta,\R}$ of the form $X_\Delta$ to the space $E_\Delta$ is real valued,
 i.e.\ $X_{\Delta,\R}\in\bigwedge_\R^{2n-k}E_\Delta^*$;

 {\rm(2)}
 $X_\Delta(\xi_1,\cdots,\xi_k)=0$ if $\exists i\colon \xi_i\in\C_\Delta$.
  \label{dfCtropPolPol}
 \end{df}
  \begin{rem}
  Condition (2) follows from condition (1).
 \end{rem}
  \begin{rem}
  In \cite{monge} we consider the exponential tropical varieties with polynomial weights.
  So in \cite{monge} for \ETP\ $X$ the frame $X_\Delta$ is an exterior form multiplied by a polynomial in the space $E_\Delta$.
  The corresponding construction of tropical geometry see in  \cite{Est}.
 \end{rem}
 \begin{df}
  Say that a real subspace $E$ of $\C^n$ is \emph{degenerate},
  if $\codimc E_\C <\codim E$,
  where $E_\C$ is the  maximal complex subspace of $E$.
  Say that the cell $\Delta$ is degenerate,
  if the space $E_\Delta$ is degenerate.
  \label{dfDegSp}
 \end{df}
 If $k=2n,2n-1$, then any $k$-dimensional subspace is nondegenerate.
 If $k<n$, then any $k$-dimensional subspace is degenerate.
 If the $k$-dimensional subspace $E$ is nondegenerate, then  $\dim E - \dim_\R E_\C = 2n-k$.
 \begin{cor}
 Let $\Delta$ be a $k$-dimensional degenerate cell of $k$-dimensional \ETP\ $X$.
 Then $X_\Delta=0$.
 \label{corNondegcell}
 \end{cor}
 \begin{df}
  The equivalence class of \ETP\
  is called exponential tropical variety (\ETV).
  If the dimension of \ETP\ is $k$,
  then (by definition) the dimension of \ETV\ is $k-n$.
  \label{dfCtropPol}
 \end{df}
   \emph{In what follows the record $X\Rightarrow\cal P$ means that \ETP\ $X$ lies in the equivalence class \ETV\ $\cal P$}.

   Let $X$ be a \ETP.
   We define the current of measure type $\bar X$ as
  \begin{equation}
  \label{barX}
 \bar X(\varphi)=\sum_{\Delta\in X,\dim\Delta=k}\int_{\Delta} X_{\Delta,\R}\wedge\varphi.
  \end{equation}
 \begin{cor}
 \ETP\ $X,Y$ are equivalent,
 if and only if $\bar X=\bar Y$.
 \label{corEqui}
 \end{cor}
  \begin{cor}
 Let $\cal P$ be a $d$-dimensional \ETV.
 Then $\bar{\cal P}$ is a current of bedegree $\(d,d\)$.
 If $d>0$,
 the the current $\bar{\cal P}$ is closed.
 \label{corClosed}
 \end{cor}
\begin{ex} \emph{Corner loci of piecewise linear functions}.
 The continuous function $h\colon\C^n\to\R$
  is called \emph{piecewise linear},
  if it is a real polynomial of degree $1$ on any $P\in\{P\}$,
  where $\{P\}$ is a finite set of convex polytopes such that $\cup_{P\in\{P\}} P=\C^n$.
 Let the support of $(2n-1)$-dimensional polyhedral set $X$  be a corner locus of a piecewise linear function $h$.
 Then any $(2n-1)$-dimensional cell $\Delta\in X$ in locally (near any internal point of $\Delta$)
 is a common face of two halfspaces $A$ and $B$.
 Let $h_A=h|_A$ and $h_B=h|_B$.
 The ordering of the pair $(A,B)$ sets
 the coorientation of $\Delta$. The standard orientation of $\C^n$ and the coorientation
of $\Delta$ together set the orientation of the cell $\Delta$. Using this orientation put
 $X_{\Delta,\R}=d^ch_A-d^ch_B$
 (remind that $d^cg(x_t) = dg(ix_t))$.
 Easy to verify that $\bar{\cal P}=dd^ch$,
 where $X\Rightarrow\cal P$.
 For any $(n-1)$-dimensional \ETV\ $\cal P$ \emph{there exists a piecewise linear function $h$ such
 that $\bar{\cal P}=dd^ch$} (Theorem \ref{thmMonge}).
  \label{exHyper}
  \end{ex}
 \subsection{Addition of ETV}
 Let ${\cal Q}_1,{\cal Q}_2$ be $(n-k)$-dimensional \ETV\
 and $X_i\Rightarrow{\cal Q}_i$.
 There is a $k$-dimensional polyhedral set $X$ such that

 (1) $\supp X=\supp X_1\cup\supp X_2$

(2) if $\Delta\cap\Xi\not=\emptyset$, where $\Delta$ and $\Xi$ are cells of  polyhedral sets $X$ and $X_i$,
 then $\Delta\cap\Xi$ is a cell of  polyhedral set $X$.

 Let $\Delta\in X$ be a $k$-dimensional cell.
 Set $X_\Delta$ equal to the sum of (one or two) frames of the cells (one or two) of polyhedral sets $X_i$
 containing $\Delta$.
 \begin{cor}
  The polyhedral set $X$ is \ETP\ and
  $\bar X=\bar X_1+\bar X_2$.
 \label{corSum}
 \end{cor}
The equivalence class of $X$ does not depend on the choice of $X_i$ (corallary \ref{corEqui}).
Now define ${\cal Q}_1+{\cal Q}_2={\cal Q}$,
where $X\Rightarrow{\cal Q}$.
 \begin{cor}
  $m$-dimensional \ETV\ form the commutative group and the map
  ${\cal P}\mapsto\overline{\cal P}$ is injection.
 \label{corMod}
 \end{cor}
Below we consider the formally defined addition of \ETV\ of all dimensions as a graded group.
By definition, the degree of $m$-dimensional \ETV\ is $n-m$.

Using the odd volume form $\omega$ of a real vector space we can define the volume of any bounded domain $U$
as $\int_U\omega$.
Indeed, this integral does not depend on the choice of orientation.
    Say that the odd volume form $\omega$ is positive (negative) if
    $\omega$ give positive (negative) volumes of bounded domains.
  \begin{df}
Let $\Delta$ be a $k$-dimensional cell of a $k$-dimensional \ETP\ $X$.
 Say that the cell $\Delta$ is \emph{positive} (\emph{negative}),
   if the direct image of the form $X_{\Delta,\R}$ on the space $E_\Delta/\C_\Delta$
   is positive (negative) volume form.
   The construction of direct image of the odd form  $X_{\Delta,\R}$ requires the coordination of the choice of orientations
   of the spaces $E_\Delta$ and $E_\Delta/\C_\Delta$.
   We do it using the orientation of $\C_\Delta$ as the standard orientation of a complex vector space.
   \label{dfPozCell}
 \end{df}
  \begin{df}
   \ETP\ with all nonnegative cells is called positive.
   If $X\Rightarrow\cal P$ and \ETP\ $X$ is positive,
   then \ETV\ $\cal P$ also is called positive.
  \label{dfPozPol}
 \end{df}
  \begin{cor}
  Any \ETV\ is the difference of two positive \ETV.
  \label{corDiff}
 \end{cor}
 Indeed,
 let
 $X\Rightarrow\cal P$.
 Let $\Delta\in X$ be a $k$-dimensional cell with nonzero frame $X_\Delta$.
  Consider
 the single-celled \ETP\ $X^\Delta$ with the cell $E_\Delta$ and it's frame $c_\Delta X_\Delta$,
 where $c_\Delta$ is such a real number,
 that the cell $E_\Delta$ is positive.
 let $X^\Delta\Rightarrow{\cal Q}^\Delta$ and
 let $|c_\Delta|$ be sufficiently large.
 Then the \ETV\ ${\cal P}+\sum_\Delta {\cal Q}^\Delta$ is positive
 and ${\cal P}=\left({\cal P}+\sum_\Delta {\cal Q}^\Delta\right) - \left(\sum_\Delta {\cal Q}^\Delta\right).$

  Say that the $k$-dimensional cells $\Delta,\Lambda$ of $k$-dimensional polyhedral set are neighbor,
 if $\dim\Delta\cap\Lambda=k-1$.
 The set of $k$-dimensional cells is called \emph{connected},
 if for any pair of cells $\Delta,\Lambda$ of this set
 there exists a sequence of cells $\Delta_1=\Delta,\Delta_2,\cdots,\Delta_m=\Lambda$ such that
  $\forall i$ the cells $\Delta_i,\Delta_{i+1}$ are neighbor.

  Let $\Xi_1,\cdots,\Xi_q$ be maximal connected subsets of cells of $k$-dimensional polyhedral set $X$.
  Let $Y_i$ be a polyhedral set formed by cells from the subset $\Xi_i$.
  Then $\supp X=\cup_{1\leq i\leq q}\supp Y_i$.
  If $i\not=j$ then $\dim(\supp Y_i\cap\supp Y_j)<k-1$.
  If $X$ is \ETP,
  then the polyhedral sets $Y_i$ with inherited chains are \ETP\ also.
  The latest statement is the direct corollary of definition \ref{dfCtropPolPol}.

  Say that \ETP\ $Y_i$ is
  \emph{the irreducible component} of \ETP\ $X$.
  If $q=1$ then \ETP\ $X$ is said to be \emph{irreducible}.
  If $X\Rightarrow\cal P$,
  then the irreducible components of \ETV\ $\cal P$ are well defined.
 \begin{cor}
  Any \ETV\ is the sum of it's irreducible components.
  Irreducible components of positive \ETV\ are positive.
  \label{corirred}
 \end{cor}
          \subsection{ETV and convex polytops in $\C^n$}
 For a formulation of theorems \ref{thmTropStructHomog}, \ref{thmTropStruct} we need some simple geometrical facts and definitions.
 \begin{df}
 The \emph{fan of cones} is a polyhedral set such that
 any cell is a cone with a zero vertex.
 If the fan is \ETP, then we call it a homogeneous \ETP.
 The corresponding ETV\ also is called homogeneous.
 \label{dfTropFan}
 \end{df}
 Let $\gamma$ be a convex bounded polytope in ${\C^n}^*$.
 \emph{The dual cone $\Delta$ of the face} $\delta$
 of polytope $\gamma$ is, by definition,
 the set of points $z\in\C^n$ such that
 $\max_{z^*\in\gamma}\re\langle z,z^*\rangle$
 is reached at any $z^*\in\delta$.
 If $\dim\delta=m$, then $\dim\Delta=2n-m$.
 Dual cones of $m$-dimensional faces of $\gamma$
 form the $(2n-m)$-dimensional fan of cones $X^{\gamma,2n-m}$.
  \begin{lm}
 Let $U$ be an open bounded domain of $m$-dimensional
 real vector space $E$.
 Then for any orientation $\alpha$ of $E$
 there exists the only multivector $p_U(\alpha)\in\bigwedge^mE$ such,
 that $\int_U\omega=\omega(p_U(\alpha))$ for any volume form $\omega$ of the space $E$.
  \label{lmLinAlg1}
 \end{lm}
It's obviously that $p_U (-\alpha)=-p_U (\alpha)$,
where $(-\alpha)$ -- is the different from $\alpha$ orientation of $E$.
I.e.\ the multivector $p_U$ is odd.
The odd multivector $p_U$ is called \emph{a volume of domain} $U$.

 Let $E_\delta$ be a tangent space of the face $\delta$.
 The set of $m$-dimensional faces of polytope $\gamma$ form a polyhedral complex.
 The function $\delta\mapsto p_\delta\in\bigwedge^mE_\Delta\subset\bigwedge^m_\R{\C^n}^*$ on the set of $m$-dimensional cells is $m$-cochain of this complex with values in $\bigwedge_\R{\C^n}^*$.
  \begin{lm}
 The $m$-cochain $\delta\mapsto p_\delta$ is a cocycle.
  \label{lmLinAlg2}
 \end{lm}
The lemma is equivalent to the  Pascal conditions for $(m+1)$-dimensional faces of polytope $\gamma$.
 (The Pascal conditions for $k$-dimensional polytope $\gamma$ is as follows: $\sum_\delta v_\delta=0$,
 where $\delta$ is $(k-1)$-dimensional face of $\gamma$ and $e_\delta$ is an external normal of the length equal to the $(k-1)$-dimensional volume of the face $\delta$).

 The volume of $(m+1)$-dimensional convex polytope equals to the sum of volumes of its $m$-dimensional faces
 multiplied by the lengths of corresponding heights and divided by $m+1$.
 Using cocycle $p_\delta$ we can write it as
\begin{equation}\label{volReal}
    p_\delta=\frac{1}{m+1}\sum_{\theta\subset\delta,\dim\theta=m}w_\theta\wedge p_\theta,
\end{equation}
where $w_\theta$ is an arbitrary point of $m$-dimensional face $\theta\subset\delta$.
\begin{lm}
Consider the complex vector space $V$ as the real vector space $E$.
There exists the only ring homomorphism $\varrho\colon\bigwedge_\R E\to\bigwedge_\C V$
 such,
 that the map $\varrho\colon E\to V$ is the identity.
  \label{lmLinAlg3}
 \end{lm}
 Let $\delta_\C$ be a minimal complex subspace of $\C^n$ containing the real subspace $E_\delta$.
 Say that the face $\delta$ of polytope $\gamma$ is degenerate
 if $\dim\delta>\dim_\C\delta_\C$.
 Any face of dimension $0,1$ is nondegenerate.
 Any face of dimension $>n$ is degenerate.
 \begin{cor}
If the face $\delta$ is degenerate,
then $\varrho(p_\delta)=0$.
  \label{corLinAlg1}
 \end{cor}
\begin{cor}
The $m$-cochain $\delta\mapsto\varrho(p_\delta)\in\bigwedge^m_\C{\C^n}^*$
is a cocycle.
If $m>n$ then this cocycle is zero.
\label{corLinAlg2}
\end{cor}
In 2.1 we use the complex variant of formula (\ref{volReal}).
\begin{equation}\label{volComplex}
    \varrho(p_\delta)=\frac{1}{m+1}\sum_{\theta\subset\delta,\dim\theta=m}w_\theta\wedge\varrho(p_\theta),
\end{equation}
where (in contrast to the formula  (\ref{volReal})) we deal with complex multivectors.

 For $m\leq n$
 we define a \emph{homogeneous $(2n-m)$-dimensional \ETP,
 corresponding to a polytope $\gamma\subset{\C^n}^*$}.
 The cells of \ETP\ are the cones of the fan $X^{\gamma,2n-m}$.
 The frame $X^{\gamma,2n-m}_\Delta$ of cone $\Delta$ is constructed as follows.

 Set $X^{\gamma,2n-m}_\Delta=0$ if the face $\delta$ is degenerate.
 Let the face $\delta$ be nondegenerate.
 We consider the multivector $(-i)^m\varrho(p_\delta)$
 as an exterior $m$-form $W^{\gamma,2n-m}_\Delta$ on $\C^n$.
 The sign of this form depends on the choice of the face $\delta$ orientation.
  To construct the frame $X^{\gamma,2n-m}_\Delta$ from $W^{\gamma,2n-m}_\Delta$ we must establish the correspondence of
 orientations of the spaces $E_\Delta$ and $E_\delta$.


 The bilinear form $\im\langle z,z^*\rangle$ give the nodegenerate
 pairing  $E_\Delta/\C_\Delta\otimes E_\delta\to\R$.
 Consider the corresponding
 symplectic form $\omega$ on the space $E_\Delta/\C_\Delta\oplus E_\delta$
 and the orientation $\chi$
 corresponding to the volume form $\omega^m$.
 Coordinate the orientations of $E_\Delta/\C_\Delta$ and $E_\delta$ so that
 together they set the orientation  $\chi$ of the space $E_\Delta/\C_\Delta\oplus E_\delta$.
 Now coordinate the orientations of the spaces $E_\Delta/\C_\Delta$ and $E_\Delta$ by choosing
 the orientation of $\C_\Delta$ as the standard orientation of a complex vector space.
 Now the orientations of the spaces $E_\delta$ and $E_\Delta$ are agreed.

 The conditions (1) and (2) of definition \ref{dfCtropPolPol} satisfied by construction.
 The closedness of the framed polyhedral set $X^{\gamma,2n-m}$
 is equivalent to the statement of Lemma \ref{corLinAlg2}.
 Thus, for $k\geq n$ the framed fan $X^{\gamma,k}$ is homogeneous \ETP.
 Corresponding homogeneous \ETV\ also denoted $X^{\gamma, k}$.
   \begin{thm}
   For any $k$-dimensional homogeneous \ETV\ $X$ exists a finite set $\{\gamma\}$ of convex polytopes in the space ${\C^n}^*$ such that
   \begin{equation}\label{structHomog}
    X=\sum_{\gamma\in\{\gamma\}} \pm X^{\gamma,k}.
  \end{equation}
 \label{thmTropStructHomog}
 \end{thm}
 In (\ref{structHomog}) and (\ref{struct}) $tX$ is,
 by definition, the \ETP\ with the same cells as \ETP\ $X$,
 and the cell's frames multiplied by $t$.
 \begin{thm}
  For any $k$-dimensional \ETV\ $X$ exists a finite set $\{\gamma\}$ of convex polytopes in the space ${\C^n}^*$ and the set of vectors $\{a_\gamma\}$ such that
  \begin{equation}\label{struct}
    X=\sum_{\gamma\in\{\gamma\}} \pm(a_\gamma+X^{\gamma,k}),
  \end{equation}
  where $a_\gamma+X^{\gamma,k}$ is a translation of \ETV\ $X^{\gamma,k}$ by the vector $a_\gamma$.
  \label{thmTropStruct}
 \end{thm}
   \subsection{Stable intersections of ETP and multiplication of ETV.}
  Say that polyhedral sets $X$ and $Y$ are transversal,
  if for any pair of cells $\Delta\in X,\Lambda\in Y$ with nonempty intersection
  the intersection of spaces $E_\Delta$ and $E_\Lambda$ is transversal.

  Let $X,Y$ be transversal polyhedral sets of dimensions $p,q$.
  Then the pairwise intersections of cells of $X$ and $Y$ form a polyhedral set $X\cap Y$.
  The dimension of the polyhedral set
  $X\cap Y$ is equal to $p+q-2n$ (if $X\cap Y\not=\emptyset$, then $p+q\geq2n$).
  For transversal \ETP\ $X,Y$ on a polyhedral set $X\cap Y$ with $\dim X\cap Y\geq n$ we can determine the structure of \ETP\ as follows.

  Let $X,Y$ be transversal framed polyhedral sets satisfying all the conditions of definition \ref{dfCtropPolPol},
 \emph{except, perhaps, the condition of closedness}.
  Let $\Delta\in X, \Lambda\in Y$ be cells of higher dimensions and let
  $\Xi=\Delta\cap\Lambda$.
  Below we define an odd form $X_\Delta\wedge Y_\Lambda$ on the cell $\Xi$.
   Then the frame of $\Xi$ is defined as $(X\cap Y)_\Xi = X_\Delta\wedge Y_\Lambda $.

  Choose the orientations $\alpha,\beta$ of cells $\Delta,\Lambda$ and let
 $X_\Delta^\alpha, Y_\Lambda^\beta$ be exterior forms
  corresponding to the chosen orientations.
  We make an odd form $X_\Delta\wedge Y_\Lambda$,
  attributing to the form $X_\Delta^\alpha\wedge Y_\Lambda^\beta$ sign
  depending on the orientation of the cell $\Xi$:
  make the form of $X_\Delta\wedge Y_\Lambda$ positive (definition \ref{dfPozCell}),
  if forms $ X_\Delta, Y_\Lambda$ are both positive or both negative.
  Otherwise, make the form $X_\Delta\wedge Y_\Lambda$ negative.
 \begin{cor}
 If the polyhedral sets $X$ and $Y$ of dimension $p$ and $q$ are transversal and closed,
 then
 the polyhedral set $(X \cap Y)$, framed as described above, is \ETP.
  \label{corProductTransv}
 \end{cor}
 Indeed,
  let $\Upsilon=\Delta\cap\Gamma$  be a $(p+q-2n-1)$-dimensional cell of the polyhedral set $X\cap Y$,
where $\Delta$ is  $p$-dimensional cell polyhedral set $X$
and $\Gamma$ is $(q-1)$-dimensional cell of polyhedral set $Y$.
Then, using the definition of the boundary of a polyhedral set (equation (\ref{bnd})),
get
 $$(\partial (X\cap Y))_\Upsilon=\sum_{Y\ni\Lambda\supset\Gamma} X_\Delta\wedge Y_\Lambda=X_\Delta\wedge\sum_{Y\ni\Lambda\supset\Gamma}  Y_\Lambda=X_\Delta\wedge(\partial Y)_\Gamma=0.$$
 Thus, the closedness of a polyhedral set $(X\cap Y)$ is a consequence of  closedness of polyhedral sets $X,Y$.
 The remaining conditions of Definition \ref {dfCtropPolPol} are obvious.
  \begin{df}
  We call ETV ${\cal P}$ and ${\cal Q}$ transversal if there exist transversal \ETP\ $X,Y$
  such that $X\Rightarrow\cal P$ and $Y\Rightarrow\cal Q$.
    \label{dfProdTransv}
 \end{df}
 Let $X\Rightarrow\cal P$ and $Y\Rightarrow\cal Q$,
 where  \ETP\ $X$ and $Y$ are transversal.
Say that ${\cal P}{\cal Q}=\cal X$,
where $X\cap Y\Rightarrow\cal X$.
The product of transversal \ETP\ is well defined.

Also, if $\dim{\cal P} + \dim{\cal Q}<n$, then (for any \ETV) set ${\cal PQ} = 0$.
 In any case for the definition of the product of \ETV\ we do the following.

 The additive group of the space $\C^n$ acts on the set of polyhedral sets by shifts.
For any fixed pair of polyhedral sets $X,Y$ shifted polyhedral sets $X,z + Y$ are transversal for almost all $z\in\C^n$.
 \begin{thm}
 If \ETP\ $X$ and $z_i+Y$ are transversal, then for $z_i\to0$ the sequence of currents $\overline{X\cap(z_i + Y)}$
  converges to the limit current $\overline {Z}$, where $Z$ is \ETP\ of dimension ($p+q-2n$).
  The limit current does not depend on the choice of the sequence $z_i$.
 If $X\Rightarrow\cal P$ and $Y\Rightarrow\cal Q$, then the equivalence class of \ETP\ $Z$ gives \ETV\ ${\cal P}{\cal Q}$ of dimension $(p-n)+(q-n)-n$ independent on the choice of \ETP\ $X,Y$.
  \label{thmProduct1}
 \end{thm}
 The following explains that if \ETP\ $X,Y$ are positive, then \ETP\ $Z$ from the statement of Theorem \ref{thmProduct1}
 is \emph{a stable intersection} of \ETP\ $X$ and $Y$.
  \begin{df}
 Let $A,B$ be subsets of finite-dimensional vector space $E$.
 A point $x\in A\cap B$ is called \emph{stable} if any neighborhood of $x$ contains the points of the set $A\cap(y+B)$ at sufficiently small shifts $y$.
 The set of stable points of intersection is called
  a stable intersection and is denoted $A\spr B$.
 \label{dfstab}
 \end{df}
 Locally any polyhedral set coincides with a fan of cones.
 Indeed, let $\Theta\in X$ and $e\in\Int\Theta$, where $\Int\Theta$ is the interior of the cell $\Theta$ (assuming that the only point of $0$-dimensional cell is internal).
 The shifted polyhedral set $(X-e)$ in a neighborhood of zero coincides with the fan of cones.
 Denote this fan by $X_\Theta$ and call it \emph{$\Theta$-localization of polyhedral set} $X$.
 The fan $X_\Theta$ is independent of the choice of internal point $e$ of the cell $\Theta$.

 Denote by $ K_\min$ the minimal cone of $K$.
 A minimal cone of a fan is a subspace.
 Cone $K_\min$ is contained in all cones of $K$.
 Minimal cone of $X_\Theta$ is the subspace $E_\Theta$.
 For any cell $\Theta$ of \ETP\ $X$ the localization $X_\Theta$ is homogeneous \ETP,
 if we equip cells of the polyhedral set $X_\Theta$ with frames inherited from $X$.
\begin{cor}
 Let $\Delta\in X,\Lambda\in Y$ and $x\in\Int\Delta\cap\Int\Lambda$.
 Then the following conditions are equivalent

 {\rm(1)} $x\in\supp X\spr\supp Y$

 {\rm(2)} $\Delta\cap\Lambda\subset\supp X\spr\supp Y$

 {\rm(3)} $(X_\Delta)_\min\cap(Y_\Lambda)_\min\subset\supp X_\Delta\spr\supp Y_\Lambda$
    \label{corStab1}
\end{cor}
\begin{cor}
 $\dim X\spr Y\leq \dim X+\dim Y-2n$
 (it is assumed that a set of negative dimension is empty).
 \label{corStab2}
\end{cor}
 Let $K,L$ be homogeneous  \ETP,
 $\dim K=k,\,\dim L=l$,
  $\supp K\spr\supp L=K_\min\cap L_\min$ and
  polyhedral sets $K,z+L$ are transversal.
  Set $V=K_\min\cap L_\min$.
  \begin{lm}
 If \ETP\ $K,L$ be positive, then $\dim V=k+l-2n$.
  \label{lmPositiveStable}
 \end{lm}
  Assume that \ETP\ $K,L$ are positive.
  Let $\{\Xi_i(z)\}$ be a set of non-empty intersections $\Delta\cap(z+\Lambda)$,
 where $\Delta\in K$ and $\Lambda\in L$ are the cells of dimension $k$ and $l$ respectively.
 Then $\{\Xi_i(z)\}$ is a set of cells of the highest dimension of  \ETP\ $K\cap(z + L)$
 (corollary \ref{corProductTransv}).
Any $\Xi_i(z)$ equals to $V+z_i$ for some $z_i\in\C^n$.
So their frames $(K\cap L)_{\Xi_i(z)}$
can be considered as exterior forms on $\C^n$
with odd dependence on the orientation of the subspace $V$.
The sum of these forms is denoted by $W_{K, L}(z)$.
  \begin{lm}
 Let $K,L$ be positive homogeneous \ETP.
 Then the form $W_{K,L}(z)$ is constant (does not depend on $z$).
  \label{lmConst}
 \end{lm}
 Let $X,Y$ be positive \ETP, $\dim X=k, \dim L=l$.
 Let $\Theta\subset\supp X\spr\supp Y$,  $\dim\Theta=k+l-2n$
 and $\Int\Theta=\Int\Delta\cap\Int\Lambda$,
 where  $\Delta\in X,\Lambda\in Y$.
 Let $K=X_\Delta$ and $L=Y_\Lambda$.
 Assign to $\Theta$ the frame $W_{K,L}(z)$ from Lemma \ref{lmConst}.
  \begin{thm}
  Let $X,Y$ be positive \ETP, $\dim X=k, \dim L=l$.
  Then all $(k+l-2n)$-dimensional cells of polyhedral set $X\spr Y$
  with above described frames form the positive \ETP\ $X\spr Y$.
  For \ETP\ $X\spr Y$ all the statements of Theorem {\rm \ref{thmProduct1}} for a polyhedral set $Z$ are true.
 \label{thmConst}
 \end{thm}
 Now we define the multiplication of \ETV\ as follows.
 Let $X_1,X_2,Y_1,Y_2$ be positive \ETP\ and let $(X_1-X_2)\Rightarrow\cal P$, $(Y_1-Y_2)\Rightarrow\cal Q$.
 Then we set $$(X_1\spr Y_1-X_1\spr Y_2-X_2\spr Y_1+X_2\spr Y_2)\Rightarrow{\cal P}{\cal Q}.$$
 It is easy to check that the product is well-defined.
   \subsection{Bergman fans of ETV}
 %
 %
  Let $\Delta\subset\R^N$ be a $k$-dimensional convex polytope.
  Let $\Delta^\infty$ be a convex polyhedral cone
  formed by the limit points of polytopes $t\Delta$ as $t\to+0$.
  This cone lies in the subspace $E_\Delta$
  (remind: the subspace $E_\Delta$ is generated by the differences of points of the polytope $\Delta$).
  If $\Delta$ is bounded, then $\Delta^\infty=0$.
  Else $\dim\Delta^\infty\not=0$.
  Any face of the cone $\Delta^\infty$ is a cone $\Lambda^\infty$,
  where $\Lambda$ is a face of $\Delta$.

  For a polyhedral set $X$, we set $S=\cup_{\Delta\in X}\Delta^\infty$.
  Let $\varphi\in S$ and let $L(\varphi)$ be a set of cells $\Delta\in X$
  such that $\varphi\in\Delta^\infty$.
  The vectors $\varphi,\psi$ of $S$ are called equivalent,
 if $L(\varphi)=L(\psi)$.
 The equivalence classes are convex cones of dimension $\leq k$, where $k=\dim X$.
 Let $S^\infty$ be a closure of the union of equivalence classes,
 represented by $k$-dimensional cones.
 If $S^\infty=\emptyset$, then set $X^\infty=0$.
 Else there exists a $k$-dimensional fan of cones $X^\infty$ such that

 (1) $\supp X^\infty=S^\infty$

 (2) any cone $K\in X^\infty$ is contained in some equivalence class of the set $S$.

 Let $X$ be a framed $k$-dimensional polyhedral set in $\C^n$ and $K\in X^\infty$, $\dim K=k$.
 Define the frame $(X^\infty)_K$ of cone $K$ as
   $\sum_{\Delta\in L(\varphi)} X_\Delta$,
  where $\varphi$ is some interior point of the cone $K$.
 \begin{thm}
  If $X$ is \ETP,
  then $X^\infty$ is a homogeneous \ETP.
  \ETV\ corresponding to \ETP\ $X^\infty$ depends on the equivalence class of \ETP\ $X$ only.
  \label{thmBergman1}
 \end{thm}
  \begin{df}
  Homogeneous \ETV\ $X^\infty$ is called a Bergman fan of \ETV\ $X$.
  \label{dfBergman}
 \end{df}
 \begin{thm}
  The map $\beta\colon X\to X^\infty$ is a homomorphism of the ring of \ETV\ to the ring of homogeneous \ETV.
  \label{thmBergman2}
 \end{thm}
  \begin{cor}
  If \ETV\ $X$ is represented as {\rm(\ref{struct})}, then
  $X^\infty=\sum_{\gamma\in\{\gamma\}} \pm X^{\gamma,k}$.
 \label{corBergman1}
 \end{cor}
  \begin{cor}
  Homomorphism $\beta$ is invariant under the action $z\colon X\mapsto(z+X)$
  of the additive group of $\C^n$ on the ring of \ETV.
 I.e.\ $(z+X)^\infty = X^\infty$ for any $z\in\C^n$.
  \label{corBergman2}
 \end{cor}
  \begin{thm}
  If $X$ is a nonzero positive \ETV,
  then also $X^\infty$ is a nonzero positive \ETV.
  \label{thmBergman3}
 \end{thm}
 \begin{cor}
 Let $X,Y$ be positive \ETV. Thus $XY\not=0$,
 if and only if $X^\infty Y^\infty\not=0$.
 \label{corBergman3}
\end{cor}
  \begin{cor}
 Let $X_1,\cdots X_m$ be positive \ETV.
 Then $X_1\cdots X_m=0$,
 if and only if $(a_1+X_1)\cdots(a_m+X_m)=0$
 for any $a_i\in\C^n$.
 \label{corBergman4}
\end{cor}
 \section{Mixed complex Monge-Ampere operator}
  \subsection{Monge-Ampere operator and ETV}
 The mixed complex Monge-Ampere operator of degree $k$ is (by definition)
 the map $(h_1,\cdots,h_k)\mapsto dd^ch_1\wedge\cdots dd^ch_k$
  (remind that $d^cg(x_{\rm tang})=dg(ix_{\rm tang})$,
  where $g$ is a real function on a complex manifold $M$).
  Below the map $(h_1,\cdots,h_k)\mapsto dd^ch_1\wedge\cdots dd^ch_k$
  is called the Monge-Ampere operator.

  If $h_1,\cdots,h_k$ are continuous convex functions on $\C^n$,
  then \cite{BT} the Monge-Ampere operator value $dd^ch_1\wedge\cdots\wedge dd^ch_k$
  is well defined as a current
  (that is a functional on the space of smooth compactly supported differential $(2n-2k)$-forms).
  This means that if the sequence of smooth convex functions $(f_i)_j$ converges locally uniformly to $h_i$,
  then the sequence of currents $dd^c(f_1)_j\wedge\cdots dd^c(f_k)_j$
  converges to the limit current,
  independent on the choice of approximation.
  This limit current is the current of measure type,
  ie it may be continued to a functional on the space of continuous compactly supported forms.
   It follows that any polynomial in the variables $dd^cg_1,\cdots,dd^cg_q$
  with continuous convex functions $g_i$ is well defined
  as a current of measure type.
  It is easy to prove that any piecewise linear function can be written as a difference of two
  convex piecewise linear functions.
  It follows that the action of Monge-Ampere operator and above defined currents are well defined for any piecewise linear functions also.
 \begin{rem}
 For piecewise linear (not necessarily convex) functions the current
 $dd^ch_1\wedge\cdots\wedge dd^ch_k$ depends only on the product of functions
 $h_1\cdots h_k$ \cite{monge}.
 This property is a very specific for piecewise linear functions.
 However, it is partially retained for piecewise pluriharmonic functions on a complex manifold \cite{Pluri}.
 \end{rem}
  Let $B$ be a ring generated by the currents $dd^cg$,
  where $g$ are piecewise linear functions on $\C^n$.
  Attach to the ring $B$ the current $\varphi\mapsto\int_{\C^n}\varphi$,
  where $\varphi$ is volume form with compact support on $\C^n$
  (the latest current is a unit of the ring $B$).
    The ring $B$ is graded by degrees of corresponding polynomials in the variables $dd^cg$.
   \begin{ut} {\rm(\cite{monge})}
   For any homogeneous element $b$ of the ring $B$ there exist $(n-\deg b)$-dimensional \ETV\ $\iota(b)$ such that $b=\overline{\iota(b)}$.
  \label{utMonge}
 \end{ut}
 \begin{thm}
 Let $\cal B$ be a graded ring of \ETV.
 Then the map $\iota\colon B\to{\cal B}$ is an isomorphism of graded rings.
  \label{thmMonge}
 \end{thm}
 \begin{rem}
 The theorem was suggested as a conjecture in \cite{monge}.
 \end{rem}
 \begin{df}
 (1)\ Function $g\colon\C^n\to\R$ is said to be positively homogeneous of degree $1$,
 if $\forall\lambda\geq0\colon g(\lambda z)=\lambda g(z)$.\\
 (2)\ Function $g\colon\C^n\to\R$ is said to be $\R$-generated,
 if $g(z+y)=g(z)$ for any $y\in\im\C^n$.\\
 (3)\ \ETV\ $X$ is said to be $\R$-generated,
 if $y+X=X$ for any $y\in\im\C^n$
 (here $y+X$ is a shift of \ETV\ by the vector $y$).
 \end{df}
 The following statements follow from Theorem \ref{thmMonge}.
   \begin{thm}
   Let ${\cal B}_\h$ be a ring of homogeneous \ETV,
   and $B_\h$ be a ring,
   generated by currents $dd^cg$,
   where $g$ are positively homogeneous of degree $1$ piecewise linear functions on $\C^n$.
 Then the map $\iota|_{B_\h}$ is an isomorphism of graded rings $\iota|_{B_\h}\colon B_\h\to{\cal B}_\h$.
 \label{thmMongeHomog}
 \end{thm}
 %
 %
 \begin{thm} {\rm(\cite{cfan})}
   Let ${\cal B}_\R$ be a ring of $\R$-generated \ETV\ and $B_\R$ be a ring,
    generated by currents $dd^cg$,
 where $g$ are $\R$-generated piecewise linear functions on $\C^n$.
 Then the map $\iota|_{B_\R}$ is an isomorphism of graded rings $\iota|_{B_\R}\colon B_\R\to{\cal B}_\R$.
  \label{thmMongeTrop}
 \end{thm}
 The ring ${\cal B}_\R$ coincides with the ring of tropical varieties in $\R^n$ \cite{cfan}.
  \begin{cor}
 The rings ${\cal B},{\cal B}_\h,{\cal B}_\R$  are generated by the elements of a first degree.
  \label{corMonge}
 \end{cor}
 The proof of Theorem \ref{thmMonge} is based on the following statement (\cite{monge}, Theorem 2 and Proposition 2).
   \begin{ut}
   Let $X$ be a $k$-dimensional \ETP\ and $X\Rightarrow\cal P$.
   Let $h$ be a piecewise linear function that is linear on any cell of polyhedral set $X$.
    Construct a new frame of degree $2n-k+1$ on any $k$-dimensional cell $\Delta\in X$ as
    $$Y_\Delta=d^cG_\Delta\wedge X_\Delta,$$
   where $G_\Delta$ is any linear function such that $(G_\Delta)|_\Delta=h$.
   Set $D_c(hX)=\partial Y$.
   Then

   {\rm(1)} $(k-1)$-dimensional framed polyhedral set $D_c(hX)$ is \ETP\ and does not depend on the choice of functions $G_\Delta$

   {\rm(2)}
   \ETV, corresponding to \ETP\ $D_c(hX)$, depends only on the equivalence class of \ETP\ $X$

   {\rm(3)} $\overline{D_c(hX)}=dd^c(h\bar{\cal P})$
  \label{utMonge2}
 \end{ut}
 It is obvious, that the map $\iota$ is bijective on the elements of degree $0$.
 Using induction on degrees of graded ring $B$ and applying statement \ref{utMonge2},
 obtain the statement \ref{utMonge} (details are in \cite {monge}).
 Thus, the map $\iota$ is defined as a homomorphism of graded abelian groups.
 By definition the map $\iota$ is injective.
 Therefore, the proof of bijectivity of $\iota$
 is reduced to the following statement.
   \begin{thm}
   For any $(k-1)$-dimensional \ETV\ $\cal P$ there exist finite sets of
   \\
(a) $k$-dimensional \ETV\ $\{\cal Q\}$
\\
(b) piecewise linear functions $\{h_{\cal Q}\}$ on the space $\C^n$
\\
(c) \ETP\ $\{X_{\cal Q}\Rightarrow\cal Q\}$
\\
  such that any function $h_{\cal Q}$ is linear on each cell of \ETP\ $X_{\cal Q}$ and
  $\sum_{\cal Q\in\{\cal Q\}} D_c(h_{\cal Q}{\cal Q})={\cal P}$.
  \label{thmdc}
 \end{thm}
  Let $h_\gamma$ be a support function of convex polytope $\gamma\subset\(\C^n)^*$
 (remind: support function is a piecewise linear function $h_\gamma(z)=\max_{z^*\in\gamma}\re\langle z,z^*\rangle$ on $\C^n$).
 From the definition follows that $\forall k$ the function $h_\gamma$ is linear on any cell of \ETP\ $X^{\gamma,k}$.
   \begin{utt}
 $D_c(h_\gamma X^{\gamma,k})=(2n-k+1)X^{\gamma,k-1}$.
  \label{uttdch}
 \end{utt}
 \proof
 Let $\Theta$ be a $k$-dimensional cone dual to a face $\theta$ of polytope $\gamma$.
 Then the frame $\left(D_c(h_\gamma X^{\gamma,k})\right)_\Delta$ of $(k-1)$-dimensional cone $\Delta$
 (by definition of the map $D_c$),
 equals to
 $$\sum_{\Theta\supset\Delta}d^cG_\Theta\wedge X^{\gamma,k}_\Theta,$$
  where $G_\Theta$ is a linear function such that $(G_\Theta)|_\Theta=h_\gamma$.
  Set $G_\Theta(z)=\re\langle z,w_\theta\rangle$,
 where $w_\theta$ is an arbitrary point of the face $\theta$.
 Then $d^cG_\Delta=\re\langle dz,-iw_\delta\rangle$.

  Let $\delta$ be a face of polytope $\gamma$ dual to cone $\Delta$.
  Then (by definition) $X^{\gamma,k-1}_\Delta=(-i)^{2n-k+1}\varrho(p_\delta)$.
  Applying (\ref{volComplex}), we get
    \begin{multline*}
 X^{\gamma,k-1}_\Delta=(-i)^{2n-k+1}\varrho(p_\delta)=\\
 \frac{1}{2n-k+1}\sum_{\theta\subset\lambda,\dim\theta=2n-k}(-iw_\theta)\wedge\left((-i)^{2n-k+1}\varrho(p_\delta)\right)=\\
 \frac{1}{2n-k+1}\sum_{\Theta\supset\Delta,\dim\Theta=k}d^cG_\Delta\wedge X^{\gamma,k}_\Theta=\frac{1}{2n-k+1}\left(D_c(h_\gamma X^{\gamma,k})\right)_\Delta
    \end{multline*}
\QED
 \textbf{Proof of Theorem \ref{thmdc}.}
 Theorem \ref{thmTropStruct} \ETV\ $\cal P$ can be represented as
     $${\cal P}=\sum_{\gamma\in\{\gamma\}} \pm(a_\gamma+X^{\gamma,k-1})$$

 Set
 ${\cal Q}=\frac{1}{2n-k+1}(a_\gamma+X^{\gamma,k})$
 (here $tX$ is,
 by definition, the \ETV\ $X$ with the frames multiplied by $t$).

 Set
 $h_{\cal Q}=h_\gamma(z-a_\gamma)$,
 where $h_\gamma$ be a support function of polytope $\gamma$.
 The function $h_{\cal Q}$ is linear on any cell of \ETP\ $X^{\gamma,k}$ by construction.
 It follows from Proposition \ref{uttdch},
 that
 $$a_\gamma+X^{\gamma,k-1}=\frac{1}{2n-k+1}D_c(h_\gamma(a_{\cal Q}+X^{\gamma,k}))=D_c(h_{\cal Q}{\cal Q}).$$
 Thus ${\cal P}=\sum_{{\cal Q}\in\{\cal Q\}} \pm D_c(h_{\cal Q}{\cal Q})$.
 \QED
 Now it remains to show that \emph{the map $\iota$ preserves the products of \ETV}
 (ie\ is an isomorphism of rings).

 The additive group of space $\C^n$ acts by shifts on the set of piecewise linear functions.
 This action extends to a continuous action on the ring $B$.
 Similarly (Theorem \ref{thmProduct1}),
 this group acts on the ring $\cal B$.
 These actions commute with the map $\iota$.
 Therefore, the statement on the isomorphism of the rings is reduced to the case of transversal intersections
 of corner loci (Proposition \ref{uttProdddc}).

 Let \ETP\ $X_1,\cdots,X_k$ be corner loci of piecewise linear functions $h_1,\cdots,h_k$.
 Say that \ETP\ $X_1,\cdots,X_k$ are transversal,
 if for all sets of cells $\Delta_i\in X_i$ with $\Delta_1\cap\cdots\cap \Delta_k\not=\emptyset$
 the tangent spaces $E_{\Delta_i}$ of cells $\Delta_i$ are transversal.
 \begin{utt}
 If \ETP\ $X_1,\cdots,X_k$ are transversal
 then $$\iota(dd^ch_1\wedge\cdots\wedge dd^ch_k)={\cal P}_1\cdots{\cal P}_k,$$
 where $X_i\Rightarrow{\cal P}_i$.
  \label{uttProdddc}
 \end{utt}
 \proof
 Let
 $\Delta_i\in X_i$ be a $(2n-1)$-dimensional cell.
 Consider a single-celled \ETP\ $A_i$ with the cell $E_{\Delta_i}$ framed as $(X_i)_\Delta$.
 Set $A_i\Rightarrow{\cal A}_i$.

 Locally $\Delta_i$ looks as a common hyperplane of two halfspaces $B_1,B_2$.
 Let a piecewise linear function $g_i$ is linear on halfspaces $B_1,B_2$ and
 $g_i|_{B_j}=h_i|_{B_j}$ near any inner point of the cell $\Delta_i$.

 It follows from transversality condition that the statement of Proposition \ref{uttProdddc} is equivalent to the series of equations
 $$\iota(dd^cg_1\wedge\cdots\wedge dd^cg_k)={\cal A}_1\cdots{\cal A}_k,$$
  where $\Delta_1\in X_1,\cdots,\Delta_k\in X_k$ is any set of $(2n-1)$-dimensional cells.
  This last equality is easy to verify directly.
  \subsection{The zero values of mixed Monge-Ampere operator at piecewise linear functions}
 \begin{df}
 Let $H$ be a complex $p$-dimensional subspace of $\C^n$
 and let $k+p>n$.
 The set of piecewise linear functions $g_1,\cdots,g_k$ on $\C^n$
 is called $H$-degenerate,
 if there exist linear functions $\varphi_i\colon\C^n\to\R$
 such that
 $\varphi_i(z+h)+g_i(z+h)=\varphi_i(z)+g_i(z)$ for any $z\in\C^n$ and any $h\in H$
 (i.e.\
 functions $\varphi_i+g_i$ are the pullbacks of some piecewise linear functions on the space $\C^n/H$ by the projection $\C^n\to\C^n/H$).
 The set $g_1,\cdots,g_k$ is called nondegenerate if there is no subspace $H$ such that the set $g_1,\cdots,g_k$ is $H$-degenerate.
 \label{dfDegLin}
 \end{df}
 \begin{thm}
  Let $h_1,\cdots,h_k$ be convex piecewise linear functions.
  Then $dd^ch_1\wedge\cdots\wedge dd^ch_k\not=0$ if and only if the set $h_1,\cdots,h_k$ is nondegenerate.
  \label{thmMongeDeg1}
 \end{thm}
 \begin{qu}
  Let $h_1,\cdots,h_k$ be any convex functions on $\C^n$.
  Is it true statement of theorem \ref{thmMongeDeg1}?
  \end{qu}
 The theorem \ref{thmMongeDeg1} is a successor of the criterion of a zero value of mixed volume.
 We explain this more.
 \begin{ut} {\rm\cite{cfan}}
 Let $A_1,\cdots,A_n$ be compact convex bodies in $\re\C^{n*}$
 and let $h_i\colon\C^n\to\R$ be their support functions.
 Consider some Euclidean metric on $\re\C^n$, the corresponding
  Hermitian metric on $\C^n$ and the dual Hermitian metric on $\C^{n*}$.
  Then the current $dd^ch_1\wedge\cdots\wedge dd^ch_n$
  is a Euclidean measure on the space $\im\C^n$ multiplied be
  the mixed volume of convex bodies $A_1,\cdots,A_n$ and by $n!$
   \label{utVolume}
 \end{ut}
The criterion of zero value of mixed volume
is formulated as follows.
 \begin{cor}
 Following statements are equivalent:

  {\rm(1)}
  the mixed volume of convex bodies $A_1,\cdots,A_N$ in $\R^n$ is zero

  {\rm(2)}
  there exist $p\leq n$ and a subset $B_1,\cdots,B_p$ of the set $A_1,\cdots,A_n$ such that
  for some $a_i\in\R^n$ all shifted bodies $a_1+B_1,\cdots,a_p+B_p$ are contained in some $(p-1)$-dimensional
  subspace of $\R^n$.
\label{corMongeDeg1}
 \end{cor}
 Let $A_i$ be a convex polytope.
 Then the corollary \ref{corMongeDeg1} follows directly from Theorem \ref{thmMongeDeg1}.
 Indeed, put the sets $A_1,\cdots,A_n$ into the space $\re{\C^n}^*$ and (using the statement \ref{utVolume}) apply
 Theorem \ref{thmMongeDeg1} to the support functions of polytopes.
  For any convex bodies $A_i$ corollary follows from \emph{the monotonicity of the mixed volume}:
 if $B_i\subseteq A_i$, then $V(B_1,\cdots,B_n)\leq V(A_1,\cdots,A_n)$.

 Here are three auxiliary statements, used in the proof of Theorem \ref{thmMongeDeg1}.
  The first is a direct consequence of the definition of \ ETV,
  second is obvious.
  Proof of the third statement is given after the proof of Theorem \ref{thmMongeDeg1}.
 \begin{lm}
 Let $X_1,\cdots,X_k$ be single-celled nonzero $(2n-1)$-dimensional \ETP and let $X_i\Rightarrow{\cal P}_i$.
 Let $\Delta_i$ be a $(2n-1)$-dimensional cells of \ETP\ $X_i$ and let $E_i$ be a tangent spaces of the cells $\Delta_i$.
 Let $0\not=\varphi_i\in{\C^n}^*$ be an equation of the hyperplane $E_i$,
 i.e.\ $\re\langle z,\varphi_i\rangle=0$ for any $z\in E_i$.
 Then $\prod_{i=1,\cdots,k}{\cal P}_i=0$,
 if and only if the vectors $\varphi_1,\cdots,\varphi_k$ are linearly dependent over $\C$.
 \label{lmcell}
\end{lm}
\begin{lm}
 Let $h$ be a piecewise linear function on $\R^N$.
 Let $v_i\in{\R^N}^*$ be differentials of function $h$ in it's areas of linearity.
 Suppose that for any $i,j$ the vectors $v_i-v_j$ are contained in some fixed subspace $H\subset{\R^N}^*$.
 Then the function $h-v_1$ is a pullback of some piecewise linear function on $\R^N/H^\bot$
 by the projection $\R^N\to\R^N/H^\bot$,
 where the subspace $H^\bot$ is the orthogonal complement of subspace $H$.
 \label{lm1cell}
\end{lm}
 Let $A_1,\cdots,A_k$ be nonempty finite sets of $n$-dimensional vector space $E$ (over arbitrary field).
 If $0<k\leq n$ and any set of vectors $a_1\in A_1,\cdots,a_k\in A_k$ is linearly dependent, then
 the \emph{set $A_1,\cdots,A_k$ is called degenerate}.
 \begin{utt}
 Let a set $A_1,\cdots,A_k$ is degenerate.
 The there exist
 \\
 (a) $1\leq p\leq k$
 \\
 (b) $(p-1)$-dimensional subspace $H\subset E$
 \\
 (c) subset $B_1,\cdots,B_p$ of the set $A_1,\cdots,A_k$
 \\
 such that $\forall j\colon B_j\subset H$.
 \label{uttfinite}
 \end{utt}
 \textbf{Proof of Theorem \ref{thmMongeDeg1}}.
 Let \ETP\ $X_i$ be a corner locus of function $h_i$.
 Let $A_i\subset{\C^n}^*$ be a set of equations of real hyperplanes $E_\Delta$
for all $\Delta$ from the set of $(2n-1)$-dimensional cells of $X_i$
(see formulation of lemma \ref{lmcell}).

From the convexity of functions $h_i$ it follows that \ETV\ $X_i$ are positive.
Therefore, applying Corollary \ref{corBergman4} and Lemma \ref{lmcell}
we get the following:
 $dd^ch_1\wedge\cdots\wedge dd^ch_k=0$,
\emph{if and only if the set of finite sets $A_1,\cdots,A_k$ is degenerate}.

If the set $A_1,\cdots,A_k$ is degenerate,
then applying Proposition \ref{uttfinite},
obtain a subset $g_1,\cdots,g_p$ of a set of functions $h_1,\cdots,h_k$ and
$(p-1)$-dimensional complex subspace $H$ of the space ${\C^n}^*$ such that
each of the functions $g_i$ satisfies the condition of Lemma \ref {lm1cell}.
Now theorem  \ref{thmMongeDeg1} follows from Lemma \ref{lm1cell}.
  \QED
%
 \textbf{Proof of Proposition \ref{uttfinite}}.
 Let $\bar C=\{C_1,\cdots,C_m\}$ be a maximal nondegenerate subset of the set $\bar A=\{A_1,\cdots,A_k\}$.
 The nondegeneracy of the set $\bar C$ implies existence of linearly independent vectors $c_1\in C_1,\cdots,c_m\in C_m$.
 Lemma \ref{lmFin1} follows from the maximality of subset $\bar C$.
 \begin{lm}
 Let $V$ be an $m$-dimensional subspace,
 generated by vectors $c_1,\cdots,c_m$,
 $C\in\bar A\setminus\bar C$, and
 $L$ be a subspace, generated by vectors of the set $C$.
 Then $L\subset V$.
  \label{lmFin1}
 \end{lm}
 \begin{lm}
 Let $V_i$ be an $(m-1)$-dimensional subspace of $E$,
  generated by vectors
$c_1,\cdots,c_{i-1},c_{i+1},\cdots,c_m$.
 Then, if $C_i\not\subset V$, then $L\subset V_i$.
  \label{lmFin2}
 \end{lm}
  \proof
 Let $\tilde V$ be a subspace generated by vectors $c_1,\cdots,\tilde{c}_i,\cdots,c_{m-1},c_m$,
 where $C_i\ni \tilde c_i\not\in V$.
 Then $\dim\tilde V=m$.
 So (Lemma \ref{lmFin1}) $L\subset\tilde V$.
 It remains to see that $V_i=V\cap\tilde V$.
 \QED
  Sort now the collection of sets $\bar C=\{C_1,\cdots,C_m\}$ so that
  $C_i\subset V$ if $i\leq l$ and
  $C_i\not\subset V$ if $i>l$.
  If $l=m$,
  then the subsets $C_1,\cdots,C_m,C$ are contained in $m$-dimensional subspace $V$
  as it required in Proposition \ref{uttfinite}.
  \begin{lm}
  Let $Q$ be an $l$-dimensional subspace
  generated by vectors $c_1,\cdots,c_l$.
  Then $L\subset Q$.
  \label{lmFin3}
  \end{lm}
  \proof
  By Lemma \ref{lmFin2}, $L\subset V_{l+1}\cap\cdots\cap V_m=Q$.
    \QED
 Note that $l>0$,
 as \ otherwise, by Lemma \ref{lmFin3},
 $L=0$, which is impossible.

 If $C_i\subset Q$ for $i=1,\cdots,l$,
 then the statement of Proposition \ref{uttfinite} is true for $p=l+1$,
 for a subset $C_1,\cdots,C_l,C$ of the set $\bar A$ and for the subspace $Q$.
 Suppose that $\exists i\leq l$ such that $C_i\not\subset Q$.
  \begin{lm}
  Let $i\leq l$ and $C_i\not\subset Q$.
  Then $L\subset Q_i$,
  where $Q_i$ is a $(l-1)$-dimensional subspace,
  generated by vectors $c_1,c_2,\cdots,c_{i-1},c_{i+1},\cdots,c_l$.
  \label{lmFin4}
  \end{lm}
  \proof
  Let $C_i\ni\breve{c}_i\not\subset Q$.
  As $i\leq l$, then $\breve{c}_i\in V$.
  Consider the decomposition
\begin{equation}
\breve{c}_i=\sum_{j\leq m} \alpha_j c_j
  \label{e3}
 \end{equation}
 respect to the basis $\{c_j\}$ of $V$.
 There may be cases $\alpha_i\not=0$ and $\alpha_i=0$.

 In the first case, the vectors
\begin{equation}
  c_1,\cdots,c_{i-1},\breve{c}_i,c_{i+1},\cdots,c_m
  \label{e4}
 \end{equation}
 are linearly independent (i.e.\ form the basis of the space $V$).
 Then, using Lemma \ref{lmFin3}, we obtain that
 $L$ belongs to the subspace generated by the vectors $c_1,\cdots,c_i,\cdots,c_{l}$
 and to the subspace generated by the vectors $c_1,\cdots,c_i,\cdots,c_{l}$.
 Hence $L$ belongs to their intersection,
 ie to the subspace $Q_i$, as required.

 In the second case (when $\alpha_i=0$) the vectors (\ref{e4}) are linearly dependent.
 In this case, since $\breve{c}_i\not\in Q$, then $\alpha_j\not=0$ for some $j>l$.
 Let $V\not\ni\tilde c_j\in C_j$.
 Then the vectors
\begin{equation}
 c_1,\cdots,c_{i-1},\breve{c}_i,c_{i+1},\cdots,c_l,\cdots,c_{j-1},\tilde c_j,c_{j+1},\cdots,c_m
  \label{e5}
 \end{equation}
 are linearly independent.
 Then (Lemma \ref{lmFin2}), $L$ belongs to subspace $\tilde Q$,
 generated by vectors
  $$c_1,\cdots,c_{i-1},\breve{c}_i,c_{i+1},\cdots,c_{j-1},c_{j+1},\cdots,c_m.$$
Then $c_i\not\in \tilde Q$.
Indeed, the otherwise we would obtain a decomposition of the form (\ref{e3}) with non-zero $\alpha_i$.
It follows that $Q\cap\tilde Q=Q_i$ and so $L\subset Q_i$.
Lemma is proved.
\QED
Now sort the set $C_1,\cdots,C_l$ so that
$C_{\leq q}\subset Q$ and $C_{>q}\not\subset Q$.
Then the space $L$ belongs to the subspace,
generated by vectors $c_{1},\cdots,c_q$.
This follows from Lemma \ref{lmFin4}, because this last subspace
coincides with the subspace $Q_{q+1}\cap\cdots\cap Q_l$.
So $q+1$ sets $C_1,\cdots,C_q,C$ are in $q$-dimensional subspace.
Proposition \ref{uttfinite} is proved.
%

\begin{flushleft}
{\small
Institute for Information Transmission Problems,
B.Karetny per. 19, 101447 Moscow, Russia\\
E-mail address: kazbori@iitp.ru
}
\end{flushleft}
\end{document}